\documentclass[12pt]{article}
\usepackage{amssymb}
\usepackage{amsthm}

\usepackage{extarrows}
\usepackage{cases}
\usepackage{graphicx}
\usepackage{mathrsfs}
\usepackage{graphics}
\usepackage[misc]{ifsym}
\usepackage{color}
\usepackage{indentfirst}
\usepackage{times}
\usepackage{cite}

\numberwithin{figure}{section}
\numberwithin{equation}{section}

\newtheorem{theorem}{Theorem}[section]
\newtheorem{proposition}[theorem]{Proposition}
\newtheorem{definition}[theorem]{Definition}
\newtheorem{corollary}[theorem]{Corollary}
\newtheorem{lemma}[theorem]{Lemma}
\newtheorem{remark}[theorem]{Remark}
\newtheorem{example}[theorem]{Example}

\newcommand{\cA}{{\mathcal A}}

\newcommand{\cC}{{\mathcal C}}

\newcommand{\cF}{{\mathcal F}}

\newcommand{\cV}{{\mathcal V}}

\newcommand{\sC}{{\mathscr C}}

\newcommand{\sG}{{\mathscr G}}

\newcommand{\trl}{{\triangleq}}

\def\N{\mathbb{N}}

\def\1{\mathbb{1}}

\def\disp{\displaystyle}

\def\bc{\begin{center}}
\def\ec{\end{center}}
\def\be{\begin{equation}}
\def\ee{\end{equation}}
\def\ba{\begin{array}}
\def\ea{\end{array}}
\def\benu{\begin{enumerate}}
\def\eenu{\end{enumerate}}
\def\bt{\begin{theorem}}
\def\et{\end{theorem}}
\def\bl{\begin{lemma}}
\def\el{\end{lemma}}
\def\bco{\begin{corollary}}
\def\eco{\end{corollary}}
\def\bn{\begin{numcases}}
\def\en{\end{numcases}}
\def\br{\begin{remark}}
\def\er{\end{remark}}
\def\bd{\begin{definition}}
\def\ed{\end{definition}}
\def\bp{\begin{proposition}}
\def\ep{\end{proposition}}
\def\bo{\begin{proof}}
\def\eo{\end{proof}}
\def\bx{\begin{example}}
\def\ex{\end{example}}

\def\De{\Delta} \def\de{\delta}

\def\sig{\sigma}
\def\vp{\varphi}

\def\~{\widetilde}

\def\ol{\overline}

\def\ra{\rightarrow}

\def\stac{\stackrel}
\def\limra{\lim_{\longrightarrow}}
\def\limla{\lim_{\longleftarrow}}
\def\8{\infty}
\def\X{\times}

\def\mb{\mbox}

\def\sm{\setminus}

\def\ss{\subset}

\def\Hs{\hspace{0.8cm}}
\def\hs{\hspace{0.4cm}}
\def\Vs{\vskip10pt}
\def\vs{\vskip5pt}

\def\[{\left[}
\def\]{\right]}
\def\({\left(}
\def\){\right)}

\begin{document}
\title{A Theory of Compact Hausdorff Shape in Hausdorff Spaces}


\author{Jintao Wang\\
{\small Center for Mathematical Sciences, Huazhong University of Science and Technology,}\\{\small Wuhan 430074, China (wangjt@hust.edu.cn)}\\
Jinqiao Duan\\
{\small Department of Applied Mathematics, Illinois Institute of Technology,}\\{\small Chicago, IL 60616, USA (duan@iit.edu)}}





\maketitle
\begin{abstract} Shape is a more general concept of homotopy type for topological spaces. This work aims to establish a new shape theory, i.e., compact Hausdorff shape (CH-shape) in general Hausdorff spaces. We use an `internal' method and direct system approach on the homotopy category of compact Hausdorff spaces. This new shape  preserves most   properties of compactly generated shape (H-shape) given by Rubin and Sanders. More importantly, it allows to develop the entire (co)homology theory for CH-shape, except for the exactness of cohomology theory, and this is dual to the approach and consequences of Marde\v si\'c and Segal.
\end{abstract}

\section{Introduction}
Shape is a more general concept of homotopy type to study the geometric properties of complicated topological spaces. The concept `shape' was originally introduced by Polish topologist Borsuk at the `International Symposium on Topology and its Applications' in  1968. Indeed in his paper `Concerning homotopy properties of compacta' \cite{Bor1}, Borsuk had actually already founded the the theory of shape of a compactum. Both Borsuk \cite{Bor2,Bor3} and Fox \cite{Fox} extended shape to arbitrary metric spaces in the following years. The method used by Borsuk seems somehow to be `internal' in nature, while Fox's method can be considered to be `external'.

Fox \cite{Fox} introduced the inverse system to define his shape equivalence, and his definition of shape is much coarser than Borsuk's. In particular, Borsuk shape can not be preserved by sums and products \cite{Bor3}. Marde\v si\'c and Segal \cite{Mar3, Dyd} developed the shape theory after Fox's method and established a more general   theory of shape in an arbitrary category by the inverse system approach. Following Fox, Marde\v si\'c and Segal adopted the ANR-systems, which was originally presented by Borsuk, to give the definition of shape of arbitrary metric spaces and compact Hausdorff spaces \cite{Mar1,Mar2}, since any such spaces can always be embed in some ANR, Here ANR stands for the absolute neighborhood retract of metric spaces that can be found in many books about algebraic topology (see, \cite{Eil,Mar3}). The shape theory developed by Marde\v si\'c and Segal in \cite{Mar3}   includes shape (homotopy) groups, shape (co)homology groups and   other related topics, just like the homotopy theory.

Based on the above-mentioned shape theory for compact Hausdorff spaces, Rubbin and Sanders \cite{Rub} extended the definition of shape to general Hausdorff spaces, via direct systems of compact subsets in the shape category of compact Hausdorff spaces and shape maps. They called such a shape the `Compactly generated shape' and their method is an `internal' one. Sanders developed H-shape theory after 1973, including shape (homotopy) groups, the sums and products of shape and Whitehead theorem (see, \cite{Sand1,Sand2,Sand3,Sand4}).

Besides Sanders' work, there was also another type of `internal' approach to describe shape --- the proximate net approach. In 1974, Felt \cite{Felt} answered partly Klee's question about the relation of $\epsilon$-continuity and shape in compact metric spaces (see \cite{Klee}). Sanjurjo \cite{Sanj1} in 1985 presented  a full description of the shape category of compacta in terms of $\epsilon$-continuity, in which he introduced the notion of a proximate net between compacta. He later gave another description \cite{Sanj2} using multi-nets. This multivalued approach was generalized by Mor\'on and Ruiz del Portal \cite{Moron} (using normal open coverings) to obtain an internal description of the shape category of paracompact spaces. By combining the original single-valued approach in \cite{Sanj1} and the techniques in \cite{Moron}, Kieboom provided a new description of the shape of paracompact spaces in \cite{Kieb},  with generalized tools --- $\cV$-continuity ($\cV$ is an open covering) and approximate net.  In the last decade, Shekutkovski and his co-authors \cite{Shek1,Shek2,Shek3,Shek4} developed intrinsic approach to shape and presented remarkable consequences using the intrinsic approach in paracompact spaces. Especially, he showed that the proximate fundamental group is one invariant of the intrinsic shape for paracompact spaces \cite{Shek3} in 2015.
\vs

Now we follow Sanders' work on the H-shape by the `internal' approach. As is known, the inverse limit can not preserve the exactness property of a given long sequence in general. As a result, the long \v Cech homology (shape homology) sequence is not necessary to be exact except for some particular cases, including the case when the space pair is compact and the coefficient group is compact or a finite dimensional vector space over a field, and some other cases \cite{Eil,Mer}.

Since Rubin and Sanders established his H-shape via direct systems in the shape category of compact Hausdorff spaces and shape maps, if one wants to construct the corresponding H-shape homology or cohomology theory, either the direct systems of \v Cech homology groups or the inverse systems of \v Cech cohomology groups are necessary to be used. In each case, exactness property can not be generally obtained for the H-shape (co)homology theory. This limitation prevents us adopting the exactness in arbitrary Hausdorff spaces and coefficient groups.
\bigskip

In this paper, we establish a new type of shape for arbitrary Hausdorff spaces, called compact Hausdorff shape (CH-shape). We develop our shape also by the `internal' method, but in the homotopy category of compact Hausdorff spaces, whose objects are compact Hausdorff spaces and whose morphisms are homotopy classes of maps between them, not the shape maps that Rubin and Sanders used. Notice each class of maps between compact Hausdorff spaces can induce a shape map between them; our type of shape is indeed a much finer shape relative to H-shape. Therefore, most of the properties possessed by H-shape also hold for CH-shape, including sums and products of CH-shape, shape (homotopy) groups and related properties. Especially, homotopy groups are CH-shape invariant, and can be regard as CH-shape groups. What is more, one can see that the CH-shape theory is fundamentally dual to the shape theory given by Marde\v si\'c and Segal \cite{Mar3} in the inverse system approach. This construction of shape is indeed a completion of shape theory.

Particularly, singular homology groups are CH-shape homology invariant. However, dual to Marde\v si\'c's shape theory, the CH-cohomology theory is not so satisfactory, since we can not avoid the utilization of inverse systems in the definition of CH-shape cohomology groups. As a result, the CH-shape cohomology groups can not necessarily meet the exactness axiom in general. But similarly, if the Hausdorff space pair is compact and the coefficient group is compact or a finite dimensional vector space over a field, the long CH-shape cohomology sequence can also be exact. These (co)homology theories can not ensured by H-shape, but they can provide much useful  information about Hausdorff spaces and have   applications in other  issues.

To start our topic,   we introduce the concepts of direct systems and direct limit for a general category as the preliminaries in the next section. We define the CH-shape in Section 3. In Section 4, we present the basic topological properties of CH-shape, including homotopy invariance and relations with other homotopy invariants. In the final section, we discuss about CH-shape (co)homology theories in details.

\section{Preliminaries}\label{s2}

In this section, we recall   some   definitions \cite{Eil, Mar2} of direct systems of an arbitrary category that will help to develop our theory.


\Vs
\subsection{Category of direct systems}
A {\em directed set} is a preordered set $A$, provided that for any $a_1$, $a_2\in A$, there is $a\in A$ such that $a_1\leq a$ and $a_2\leq a$. A directed set $A'$ is a {\em subset} of $A$, if $a\in A'$ implies $a\in A'$ and $a_1\leq a_2$ in $A'$ implies $a_1\leq a_2$ in $A$. Let $A'$ be a subset of $(A,\leq)$. If for every $a\in A$, there exists $a'\in A'$ such that $a\leq a'$, $A'$ is said to be {\em cofinal} in $A$.

Let $\sC$ be an arbitrary category. A {\em direct system} in the category $\sC$ consists of a directed set $A$, called the {\em index set}, of an object $X_a$ from $\sC$ for each $a\in A$ and of a morphism $p_{aa'}:X_a\ra X_{a'}$ from $\sC$ for each pair $a\leq a'$ such that
\benu\item[(i)]$p_{aa}=1_a=1_{X_a}:X_a\ra X_a$,\item[(ii)]$a\leq a'\leq a''$ implies $p_{aa'}p_{a'a''}=p_{aa''}$.\eenu
We denote a direct system by $X^*=\{X_a,p_{aa'},A\}$.

A {\em morphism of direct systems} $F=(f_a,\,f):X^*\ra Y^*=\{Y_b,q_{bb'},B\}$ consists of an order-preserving function $f:A\ra B$, i.e., $a\leq a'$ implies $f(a)\leq f(a')$ in $B$, and morphisms $f_a:X_a\ra Y_{f(a)}$ from $\sC$ such that whenever $a\leq a'$,
\be\label{1.1}q_{f(a)f(a')}f_a=f_{a'}p_{aa'}.\ee
One can define an identity morphism of direct systems $1_{X^*}:X^*\ra X^*$ by considering the identity function $1_A:A\ra A$ and the identity morphisms $1_a=1_{X_a}:X_a\ra X_a$, which well satisfies (\ref{1.1}). It is also clear that the composition of two morphisms of direct systems is also a morphism of direct systems, $F1_{X^*}=F$ and $1_{Y^*}G=G$.

Thus we obtain a new category from the category $\sC$, denoted by {\em dir-$\sC$}, whose objects are all the direct systems in $\sC$ and whose morphisms are the morphisms of direct systems.

If $A'$ is a subset of $A$, and $X^*=(X_a,p_{aa'},A)$ is a direct system, then the {\em direct subsystem} ${X'}^*=(X_a,p_{aa'},A')$ over $A'$ is a direct system formed by the sets and maps of $X_a$ and $f_{aa'}$ which correspond to elements and relations in $A'$. The inclusion map $i:A'\ra A$ and identity maps $1_a:X_a\ra X_a$ form a morphism $I=(1_a,i)$ from ${X'}^*$ to $X^*$. This morphism is called the {\em injection} of the subsystem into the system.


Now we define an {\em equivalence relation} $\sim$ between morphisms of direct systems. We say $(f_a,\,f)\sim(g_a,\,g)$ provided that each $a\in A$ admits a $b\in B$ with $f(a),\,g(a)\leq b$, such that
$$q_{f(a)b}f_a=q_{g(a)b}g_a.$$
Clearly, $\sim$ satisfies reflexivity, symmetry and transitivity, and hence is an equivalence relation. Thus from the category dir-$\sC$, one obtains another new category, denoted by {\em Dir-$\sC$}, whose objects are the same with dir-$\sC$ and whose morphisms are the equivalence classes of morphisms of dir-$\sC$.

For $X^*$ and $Y^*$ in Dir-$\sC$, if there are $F:X^*\ra Y^*$ and $G:Y^*\ra X^*$ in dir-$\sC$ such that $[FG]=[1_{Y^*}]$ and $[GF]=[1_{X^*}]$, we say $X^*$ and $Y^*$ are {\em equivalent} in Dir-$\sC$, denoted $X^*\sim Y^*$.

\subsection{Direct limits}

For an arbitrary direct system $X^*=\{X_a,p_{aa'},A\}$ in dir-$\sC$, the {\em direct limit} of $X^*$ consists of an object $X^\8$ in $\sC$ and homomorphisms $p_a:X_a\ra X^\8$ such that
\be\label{2.2}p_{a'}p_{aa'}=p_a,\hs a\leq a',\ee
where the homomorphism $p_a$ is often called the {\em canonical projection}.  Moreover, if $p'_a:X_a\ra Y$ is another collection of homomorphisms with property (\ref{2.2}), then there is a unique homomorphism $g:X^\8\ra Y$ such that (see, \cite{Mar2})
$$gp_a=p'_a,\hs a\in A. $$
We denote $\disp X^\8=\limra X^*$. Clearly, the direct limit $X^\8$ of $X^*$ is unique up to a natural isomorphism.

\br\label{r2.1} It is stated in \cite{Eil,Mar2} that, if $\sC$ is the category of groups and the homomorphisms between groups, then one can always construct a direct limit of a given direct system in $\sC$. In particular, when $\sC$ is the category of abelian groups, the conclusion holds true, too.\er

Now we consider the category of groups $\sG$ and the category dir-$\sG$ of direct systems in $\sG$. Consider the morphism $F=(f_a,f):X^*\ra Y^*$, we can also define the {\em direct limit} $f^\8$ of $F$ as a morphism in $\sG$ from $X^\8$ to $Y^\8$, denoted as
$$f^\8=\limra F$$
such that for any $a\in A$, $f^\8 p_a=q_{f(a)}f_a$. By a simple observation of the direct limit of direct systems, one can obtain the uniqueness of direct limit of $F$.

Referring to \cite{Eil}, with a slight extension, we have the following conclusions:\benu\item[(i)] $F\sim G$ implies $\disp\limra F=\limra G$.
\item[(ii)] Let $F:X^*\ra Y^*$ and $G:Y^*\ra Z^*$. Then $\disp\limra(GF)=\limra G\limra F$.
\item[(iii)] $\disp\limra 1_{X^*}=1_{X^\8}$.
\item[(iv)] If $X^*\sim Y^*$, then $X^\8$ and $Y^\8$ are equivalent in $\sG$.\eenu
By this observation, $\disp\limra$ is indeed a covariant functor from Dir-$\sG$ to $\sG$.

Concerning the cofinality, we have the following consequence.
\bt\label{t2.1}\cite{Eil} Let $A'$ be cofinal in $A$, $X^*$ a direct system over $A$, ${X'}^*$ the subsystem over $A'$ of $X^*$ and $I:{X'}^*\ra X^*$ the inclusion morphism. Then $\disp i^\8=\limra I$ is an isomorphism from ${X'}^\8$ to $X^\8$.\et

All the conclusions about direct systems and direct limits above also hold for inverse systems and inverse limits  \cite{Eil,Mar2}.

\section{Compact Hausdorff Shape}

In order to introduce the compact Hausdorff shape of an arbitrary Hausdorff space, like the process of the definition of compactly generated shape  by  Rubin and Sanders \cite{Rub}, we first introduce the definition of a compact Hausdorff system.
We adopt the notations of Marde\v si\'c and Segal in \cite{Mar2}. The reader is supposed to be familiar with homotopy and homology theory.
\Vs

To define the compact Hausdorff shape, we need to present the categories of direct systems of compact Hausdorff spaces. By a map $\vp:X\ra Y$ between Hausdorff spaces $X$ and $Y$, we also mean $\vp$ is continuous.

Let HCpt be the category of compact Hausdorff spaces, whose objects are compact Hausdorff spaces and whose morphisms are homotopy classes of maps between compact Hausdorff spaces.
Then by the process given in Section \ref{s2}, we obtain the category, dir-HCpt, of direct systems in HCpt. We call a direct system in HCpt a {\em HCpt-system} and any morphism of HCpt-systems a {\em HCpt-morphism}. We always use $[\vp]$ to denote the homotopy class of a map $\vp$, and call such a homotopy class a {\em homotopy map}.

In the category dir-HCpt, the equivalence relation of two HCpt-morphisms $F=([f_a],f)$, $G=([g_a],g):\,X^*\ra Y^*$ is often said to be {\em homotopy equivalence}, denoted especially by $F\simeq G$. When $F\simeq G:X^*\ra Y^*$, we say $F$ is {\em homotopic to} $G$ from $X^*$ to $Y^*$, or $F$ and $G$ are {\em homotopic}.

Thus we can induce the category Dir-HCpt, whose objects are the same as dir-HCpt and whose morphisms are homotopy equivalence classes of HCpt-morphisms. If two direct systems $X^*$ and $Y^*$ are equivalent in Dir-HCpt, they are said to have the same {\em homotopy type}, denoted by $X^*\simeq Y^*$. The equivalence class of systems $X^*$ in Dir-HCpt is called the {\em homotopy class} of $X^*$, denoted by $[X^*]$.
\Vs

Let $X$ be a Hausdorff space. Consider the family $c(X)$ of all compact subsets of $X$ ordered by inclusion, i.e., $K\leq K'$ whenever $K\ss K'\ss X$. Let
$$C(X)=\{K,\,[i_{KK'}],\,c(X)\}$$
be a system in dir-HCpt such that $K\in c(X)$ and if $K\ss K'$, then $i_{KK'}$ is the inclusion map from $K$ into $K'$. We say $C(X)$ is the {\em HCpt-system associated with} $X$.

Now we may present the main definition as follows.

\bd\label{d3.1} (Shape)  Let $X$ and $Y$ be Hausdorff spaces, and $C(X)$ and $C(Y)$ the corresponding associated systems. Let $F:C(X)\ra C(Y)$ and $G:C(Y)\ra C(X)$ be HCpt-morphisms. If $GF\simeq1_{C(X)}$, we say $X$ is {\em shape dominated} by $Y$. Moreover, when $FG\simeq1_{C(Y)}$, we say $X$ and $Y$ have the same {\em(compact Hausdorff) shape}, written by ${\rm Sh}(X)={\rm Sh}(Y)$, and $F$ or $G$ is called a {\em shape equivalence}, where $1_{X^*}:\,X^*\ra X^*$ is the identity morphism in dir-HCpt.\ed

For the readers' convenience, we recall the ANR-shape (\cite{Mar3}) and H-shape (\cite{Rub}) below. The reader can omit it if he or she is familiar with these definitions.

ANR-shape is defined for a metric space $X$, which can be always embedded into some ANR (absolute neighborhood retract) $P$ as a subspace. The open neighborhoods of $X$ in $P$ construct an inverse system $X^*$ with the ordering $U\leq V$ and the homotopy class of the inclusion map $V\ra U$ as the morphism $i_{UV}$ whenever $V\ss U$. Two metric spaces $X$ and $Y$ have the same {\em ANR-shape} provided that the inverse systems $X^*$ and $Y^*$ associated with $X$ and $Y$, respectively, are homotopy equivalent in the category of inverse systems of ANRs. An equivalence class of the morphisms $X^*$ to $Y^*$ is called a {\em shape map} from $X$ to $Y$.
H-shape is defined for a Hausdorff space $X$. The definition is almost the same with our shape defined in Definition \ref{d3.1}, except the homotopy maps used in this paper replaced by shape maps presented above.

To distinguish with the ANR-shape (denoted by ${\rm Sh}_{\rm ANR}$) and H-shape (denoted by ${\rm Sh}_{\rm H}$), we call such a new type of shape in Definition \ref{d3.1} to be {\em compact Hausdorff shape} ({\em CH-shape} for short), denoted by ${\rm Sh}_{\rm CH}$. We denote the CH-shape simply by ${\rm Sh}$ if there is no confusion. In the following sections we will define the corresponding shape of pointed Hausdorff spaces and Hausdorff space pairs via similar processes; in each case we will only use ${\rm Sh}_{\rm CH}(\cdot)$ or ${\rm Sh}(\cdot)$ to denote the corresponding shape and their differences depend on the space type in the parentheses.

Note that we have in fact determined a new shape category whose objects are Hausdorff spaces and $\mb{Mor}(X,Y)$ is the set of all homotopy classes of morphisms from $C(X)$ to $C(Y)$ in dir-HCpt. We call such a shape category the {\em CH-shape category} and denote it by CH-Sh.

\br In topological spaces, we have topological sum and product of spaces. As for our CH-shape, we can also define the sum and product of CH-shapes of Hausdorff spaces $X^a$, $a\in A$ is the index, such that
$$\sum_{a\in A}{\rm Sh}(X^a)={\rm Sh}\(\sum_{a\in A}X^a\)$$
$$\mb{and}\hs \prod_{a\in A}{\rm Sh}(X^a)={\rm Sh}\(\prod_{a\in A} X^a\).$$
The related definitions are  similar to Rubin and Sanders' work in \cite{Rub}, and hence we omit the details.\er

\section{Basic Properties of CH-shape}

\subsection{Homotopy invariance}

We first consider the homotopy properties of CH-shape. Let $X$, $Y$ and $Z$ be Hausdorff spaces. We denote the unit interval $[0,1]$ by $I$. Similarly to the results of H-shape in \cite{Rub}, we have the following theorems.

\bt\label{t4.1} A map $\vp:X\ra Y$ induces a morphism $F:C(X)\ra C(Y)=\{L,[j_{LL'}],c(Y)\}$ in dir-HCpt, such that if $\vp,\,\psi:X\ra Y$ are homotopic maps, denoted by $\vp\simeq\psi$, then $F\simeq G:C(X)\ra C(Y)$, where $G$ is induced by $\psi$.\et

\bo Since $\vp$ is continuous, for any $K\in c(X)$, $\vp(K)\in c(Y)$. We define $f:\,c(X)\ra c(Y)$ and $f_{K}:K\ra f(K)$ for each $K\in c(X)$ such that $f(K)=\vp(K)$ and $f_K=\vp|_{K}$. Then for any $K$, $K'\in c(X)$ with $K\ss K'$, we have $f(K)\ss f(K')$ and
$$[j_{f(K)f(K')}]\circ[f_{K}]=[j_{\vp(K)\vp(K')}\circ\vp|_{K}]=[\vp|_{K'}\circ i_{KK'}]=[f_{K'}]\circ[i_{KK'}],$$
which indicates that $F\trl([f_{K}],f): C(X)\ra C(Y)$ is a morphism in dir-HCpt.

For $\phi\simeq\psi:X\ra Y$,  we denote $G=([g_k],g)$ as the induced morphism by $\psi$. Let $H:\,X\X I\ra Y$ be a homotopy between $\vp$ and $\psi$, i.e., $H(\cdot,0)=f$ and $H(\cdot,1)=g$. For an arbitrary $K\in c(X)$, defining $L=H(K,I)$, one has $L\in c(Y)$. Obviously, $f(K)$, $g(K)\ss L$ and
$$j_{\vp(K)L}\circ\vp|_K\simeq j_{\psi(K)L}\circ\psi|_K:\,K\ra L.$$
Then
$$[j_{f(K)L}]\circ[f_K]=[j_{\vp(K)L}\circ\vp|_K]=[j_{\psi(K)L}\circ\psi|_K]=[j_{g(K)L}]\circ[g_K].$$
This means $F\simeq G:C(X)\ra C(Y)$ and completes the proof.\eo

\bt\label{t4.2} Let $C(X)$, $C(Y)$ and $C(Z)$ be the systems associated with $X$, $Y$ and $Z$, respectively. If $F:C(X)\ra C(Y)$, $G:C(Y)\ra C(Z)$ are morphisms induced by $\vp:X\ra Y$ and $\psi:Y\ra Z$, respectively, then $GF:C(X)\ra X(Z)$ is induced by $\psi\vp:X\ra Z$. The identity morphism $1_{C(X)}:C(X)\ra C(X)$ is induced by the identity $1_X:X\ra X$.\et

Theorem \ref{t4.1} and \ref{t4.2} imply that the shape is a homotopy invariant.

\bt  (Homotopy invariance)   Let $X$ and $Y$ be two Hausdorff spaces and $X\simeq Y$. Then ${\rm Sh}(X)={\rm Sh}(Y)$.\et

\br The CH-shape ${\rm Sh}_{\rm CH}$ is indeed a covariant functor from HCpt to CH-Sh.\er

\subsection{Relation with other homotopy invariants}

Generally speaking, CH-shape differs from the other types of shape such as ANR-shape (\cite{Mar3}), H-shape (\cite{Rub}) and the intrinsic shape by proximate approach (\cite{Shek4}). But by the replacement of shape maps (\cite{Mar3}) in \cite{Rub} into homotopy maps in this present  paper, we have the following implication.

\bt\label{t4.3} Let $X$ and $Y$ be Hausdorff spaces. If ${\rm Sh}_{\rm CH}(X)={\rm Sh}_{\rm CH}(Y)$, then ${\rm Sh}_{\rm H}(X)={\rm Sh}_{\rm H}(Y)$.\et
\bo Since ${\rm Sh}_{\rm CH}(X)={\rm Sh}_{\rm CH}(Y)$, we have $F=([f_K],f):C(X)\ra C(Y)$ and $G=([g_L],g):C(Y)\ra C(X)$ such that
\be\label{4.1}GF\simeq 1_{C(X)}\Hs\mb{and}\Hs FG\simeq 1_{C(Y)}.\ee

By the shape theory in \cite{Mar3}, ${\rm Sh}_{\rm ANR}$ is actually a covariant functor from HCpt to (ANR-)Sh(Cpt), where (ANR-)Sh(Cpt) is the ANR-shape category of compact Hausdorff spaces.

Thus we obtain two CS-morphism $F'=(\ol{f}_K,f)$ and $G'=(\ol{g}_L,g)$ between two CS-systems (see, \cite{Rub}) $\cC(X)=\{K,\ol{i}_{KK'},c(X)\}$ and $\cC(Y)=\{L,\ol{j}_{LL'},c(Y)\}$, where we use $\ol{\vp}:={\rm Sh}_{\rm ANR}([\vp])$ to denote the shape map induced by the homotopy map $[\vp]$. Concerning the definition of homotopy equivalence for HCpt-systems, we perform ${\rm Sh}_{\rm ANR}$ on (\ref{4.1}), and then we immediately obtain
$$G'F'\simeq 1_{\cC(X)}\Hs\mb{and}\Hs F'G'\simeq 1_{\cC(Y)}.$$
The proof is finished.\eo

Now we consider the relation of CH-shape and homotopy type. For this we need to recall the concept of `CS-cofinality' (see, \cite{Rub}), a specific version of cofinality in the coverings of Hausdorff spaces.
A covering $\cF$ of a Hausdorff space $X$ is said to be {\em CS-cofinal} if there is a function $g:c(X)\ra\cF$ such that (1) if $K\in c(X)$, then $K\ss g(K)$, (2) if $K,\,K'\in c(X)$ and $K\ss K'$, then $g(K)\ss g(K')$.

For $\cF$  being a compact CS-cofinal covering of $X$, we denote $F^*=\{F,[i_{FF'}],\cF\}$ be the HCpt-system, where $\cF$ is directed by inclusion and if $F\ss F'$, then $i_{FF'}:F\ra F'$ is the inclusion map. It is  clear  that $F^*\simeq C(X)$. If additionally  $X$ is compact, we can consider only the special HCpt-system $X^*=\{X,[1_X]\}$, where the index set is a singleton. Then $\{X\}$ itself as a covering of $X$, is CS-cofinal and $X^*\simeq C(X)$.

The following theorem implies that, in a sense,  CH-shape is stronger than H-shape but weaker than homotopy type.

\bt\label{t4.4} If $X$ and $Y$ are compact Hausdorff spaces, then $X\simeq Y$ if and only if ${\rm Sh}_{\rm CH}(X)={\rm Sh}_{\rm CH}(Y)$.\et

\bo It is sufficient to show that ${\rm Sh}_{\rm CH}(X)={\rm Sh}_{\rm CH}(Y)$ implies $X\simeq Y$. Since both $X$ and $Y$ are compact, we can pick two special HCpt-systems $X^*=\{X,[1_X]\}$ and $Y^*=\{Y,[1_Y]\}$ corresponding to $X$ and $Y$, respectively. Then $X^*\simeq C(X)$ and $Y^*\simeq C(Y)$. Let $F:C(X)\ra C(Y)$ be the CH-shape equivalence. By taking the composition of CH-shape equivalences as follows
$$X^*\ra C(X)\stac{F}{\ra}C(Y)\ra Y^*,$$
we have a CH-shape equivalence $H:X^*\ra Y^*$. One sees $H$ is indeed a homotopy map $[h]:X\ra Y$ for some map $h:X\ra Y$. Associated with $H$ there is a CH-shape equivalence $H':Y^*\ra X^*$, which also possesses only one homotopy map $[h']:Y\ra X$ for some $h'$. Moreover, $[hh']=[1_Y]$ and $[h'h]=[1_X]$, which proves $X\simeq Y$.\eo

\br Since H-shape and ANR-shape are equivalent for some compact Hausdorff spaces (see, \cite{Rub}), one can easily deduce from Theorem \ref{t4.4} that ANR-shape is CH-shape invariant for compact Hausdorff spaces.\er
\subsection{Homotopy groups are CH-shape invariant}\label{s4.3}

Another important homotopy invariant is homotopy groups. Actually, homotopy groups are the CH-shape invariants. To explain this, we need to consider the pointed Hausdorff spaces and the homotopy category of pointed compact Hausdorff spaces, denoted by HCpt$_*$, whose objects are pointed compact Hausdorff spaces and whose morphisms are homotopy maps between pointed Hausdorff spaces. Then we can similarly define HCpt$_*$-systems, HCpt$_*$-morphisms, dir-HCpt$_*$, homotopy equivalences and Dir-HCpt$_*$.

As for an arbitrary pointed Hausdorff space $(X,x)$, the collection of all the compact subsets of $(X,x)$ is denoted   by $c(X,x)$ and the HCpt$_*$-system associated with $(X,x)$ is
$$C(X,\,x)=\{(K,\,x),\,[i_{KK'}],\,c(X,\,x)\},$$
where if $K\ss K'$, then $i_{KK'}:(K,x)\ra(K',x)$ is the inclusion map. In a similar manner, one can determine that two pointed Hausdorff spaces $(X,\,x)$ and $(Y,\,y)$ have the same CH-shape, denoted by ${\rm Sh}_{\rm CH}(X,x)={\rm Sh}_{\rm CH}(Y,y)$, if $[C(X,x)]=[C(Y,y)]$ in Dir-HCpt$_*$.

As is know to all, given a pointed Hausdorff space $(X,x)$, one has a homotopy group $\pi_n(X,x)$ for each $n\in\N$. Indeed the operation $\pi_n$ can be seen as a covariant functor from HCpt$_*$ to $\sG$, since for the homotopy groups, $f\simeq g:(X,x)\ra(Y,y)$ implies $f_n=g_n:\pi_n(X,x)\ra\pi_n(Y,y)$. Thus each HCpt$_*$-system $X^*=\{(X_a,x_a),[p_{aa'}],A\}$ induces for each $n\in\N$ a direct system in $\sG$
$$\pi_n(X^*)=\{\pi_n(X_a,x_a),p_{aa'n},A\},$$
where if $a\leq a'$, then $p_{aa'n}:\pi_n(X_a,x_a)\ra\pi_n(X_{a'},x_{a'})$ is the homomorphism induced by the map $p_{aa'}$,

If $F=([f_a],f):X^*\ra Y^*$ is an HCpt$_*$-morphism,  then for each $n\in\N$, $F$ induces a morphism of direct systems $\pi_n(F)=(f_{an},f):\pi_n(X^*)\ra\pi_n(Y^*)$ such that $f:A\ra B$ is an order-preserving function and if $a\in A$ then $f_{an}:\pi_n(X_a,x_a)\ra\pi_n(Y_{f(a)},y_{f(a)})$ is the homomorphism induced by $f_a:(X_a,x_a)\ra(Y_{f(a)},y_{f(a)})$. The identity HCpt$_*$-morphism $1_{X^*}:X^*\ra X^*$ induces the identity morphism of direct systems: $\pi_n(1_{X^*})=1_{\pi_n(X^*)}:\pi_n(X^*)\ra\pi_n(X^*)$ and the morphism induced by a composition is the composition of the induced morphisms. By these observations, $\pi_n$ can be generalized to be a covariant functor from the category dir-HCpt$_*$ to dir-$\sG$ for each $n\in\N$.

We also see that the functor $\pi_n$ preserves homotopy equivalences of morphisms and objects in dir-HCpt$_*$, (for similar results, see \cite{Sand1}) i.e., for each $n\in\N$,
\benu\item[(1)] if $F\simeq G:X^*\ra Y^*$, then $\pi_n(F)\simeq\pi_n(G):\pi_n(X^*)\ra\pi_n(Y^*)$;\item[(2)]if $X^*\simeq Y^*$, then $\pi_n(X^*)\simeq\pi_n(Y^*)$.\eenu
\Vs

With these preparations, we have the following consequence:
\bt\label{t4.8} Let $(X,x)$ and $(Y,y)$ be two Hausdorff spaces. If ${\rm Sh}_{\rm CH}(X,x)={\rm Sh}_{\rm CH}(Y,y)$, then $\pi_n(X,x)\approx\pi_n(Y,y)$ for any $n\in\N$.\et
\bo Corresponding to the systems $C(X,x)$ and $C(Y,y)=\{(L,y),[i_{LL'}],c(Y,y)\}$, one has two direct systems of homotopy groups:
$$\pi_n(C(X,x))=(\pi_n(K,x),i_{KK'n},c(X,x))$$$$\mb{and}\hs \pi_n(C(Y,y))=(\pi_n(L,y),j_{LL'n},c(Y,y)),$$
where $i_{KK'n}:\pi_n(K,x)\ra\pi_n(K',x)$ and $j_{LL'n}:\pi_n(L,y)\ra\pi_n(L',y)$ are induced by $i_{KK'}$ and $j_{LL'}$, respectively.
Let $F=([f_K],f):C(X,x)\ra C(Y,y)$ be a CH-shape equivalence. We have also $\pi_n(F)=(f_{Kn},f)$ such that $f_{Kn}:\pi_n(K,x)\ra\pi_n(f(K),y)$ (Note $f_K(x)=y$ for all $(K,x)\in c(X,x)$) is induced by $f_K$.

By the properties of CH-shape, functors $\pi_n$ and $\disp\limra$, it suffices to show
$$\limra\pi_n(C(X,x))=\pi_n(X,x).$$
Indeed, for $(K,x)\in c(X,x)$, let $i_{Kn}:\pi_n(K,x)\ra\pi_n(X,x)$ be the homomorphism induced by the inclusion $i_K:(K,x)\ra(X,x)$. Then for any $(K,x)\ss(K',x)$,
$$i_{Kn}=i_{K'n}i_{KK'n},$$
which is induced by the equation of inclusions $i_{K}=i_{K'}i_{KK'}$.

Now, we show for any group $G$ with homomorphisms $p_{Kn}:\pi_n(K,x)\ra G$ such that $p_{Kn}=p_{K'n}i_{KK'n}$ for $(K,x)\ss(K',x)$, there is a unique homomorphism $\sig:\pi_n(X,x)\ra G$ such that
\be\label{4.2}\sig i_{Kn}=p_{Kn}.\ee
Suppose $[f]\in\pi_n(X,x)$ has a representative $f:(S^n,1)\ra(X,x)$. Let $K=f(S^n)$ and $f':(S^n,1)\ra(K,x)$ such that $f'(y)=f(y)$ for any $y\in S^n$. Then $(K,x)\in c(X,x)$ and $[f']\in\pi_n(K,x)$. Now we set
$$\sig[f]=p_{Kn}[f'].$$

We first verify $\sig$ is a well-defined. Indeed, if $f$, $g:(S^n,1)\ra(X,x)$ satisfy $[f]=[g]$, let $h:(S^n,1)\X I\ra (X,x)$ be a homotopy between $f$ and $g$, $K'=g(S^n)$, $K''=h(S^n\X I)$ and $g':(S^n,1)\ra(K',x)$ defined by $g$. One has $K\cup K'\ss K''$ and $i_{KK''}f'\simeq i_{K'K''}g':(S^n,1)\ra(K'',x)$. Then,
$$\ba{rl}\sig[f]&=p_{Kn}[f']=p_{K''n}i_{KK''n}[f']\\&=p_{K''n}[i_{KK''}f']=p_{K''n}[i_{K'K''}g']\\&=p_{K''n}i_{K'K''n}[g']=p_{K'n}[g']\\&=\sig[g].\ea$$
Every $[f]\in\pi_n(K,x)$ for some $(K,x)\in c(X,x)$ has a representative $f:(S^n,1)\ra(K,x)$. Then $[i_{K}f]\in\pi_n(X,x)$ and $f$ is defined by $i_{K}f$. Thus we have
$$\sig i_{Kn}[f]=\sig[i_Kf]=p_{Kn}[f],$$
which confirms (\ref{4.2}).

Then we show $\sig$ is a homomorphism. Let $f$, $g:(S^n,1)\ra(X,x)$, $K=f(S^n)$, $K'=g(S^n)$ and $K''=f*g(S^n)$, where $*$ is the multiplication. Denote $f'$ and $g'$ as before defined by $f$ and $g$ respectively. Then $i_{KK''}f'*i_{K'K''}g':(S^n,1)\ra(K'',x)$ is defined by $f*g$. Thus
$$\ba{rl}\sig[f*g]&=p_{K''n}[i_{KK''}f'*i_{K'K''}g']\\&=p_{K''n}[i_{KK''}f']*p_{K''n}[i_{K'K''}g']\\
&=p_{K''n}i_{KK''n}[f']*p_{K''n}i_{K'K''n}[g']\\
&=p_{Kn}[f']*p_{K'n}[g']\\
&=\sig[f]*\sig[g].\ea$$

Finally, $\sig$ is unique. Suppose there is another homomorphism $\tau:\pi_n(X,x)\ra G$ satisfying (\ref{4.2}) with $\sig$ replaced by $\tau$. Then for any $[f]\in\pi_n(X,x)$, let $f$ be its representative and $K=f(S^n)$.
Set $f':(S^n,1)\ra(K,x)$ as defined by $f$. Then
$$\tau[f]=\tau[i_Kf']=\tau i_{Kn}[f']=p_{Kn}[f']=\sig i_{Kn}[f']=\sig[f],$$
which guarantees $\tau=\sig$ and completes the proof.\eo

\section{CH-Shape (Co)homology Groups}

In this section we  develop the (co)homology theory for the CH-shape. And this is a significant property of CH-shape better than H-shape described in Section 1.

Necessarily we consider a Hausdorff space pair $(X,X_0)$, i.e., $X_0$ is a subset of Hausdorff space $X$ and $X_0\ss X$, and maps between Hausdorff space pairs $f:(X,X_0)\ra(Y,Y_0)$, which is defined as a map from $X$ to $Y$ such that $f(X_0)\ss Y_0$. Two continuous maps $f$, $g:(X,X_0)\ra(Y,Y_0)$ are said to be homotopy equivalent, if $f\simeq g:X\ra Y$ and the image of $X_0$ is contained in $Y_0$ all the time along the homotopy. A compact Hausdorff space pair $(X,X_0)$ means that both $X$ and $X_0$ are compact.

By setting the homotopy identity $[1_{(X,X_0)}]$ and composition in the usual way, we can construct a homotopy category of compact Hausdorff space pairs, denoted by HCpt$^2$, whose objects are compact Hausdorff space pairs and whose morphisms are the homotopy maps between them. Similarly to the single space case, the definitions of HCpt$^2$-system, HCpt$^2$-morphism and the category dir-HCpt$^2$ are naturally given, as well as the homotopy equivalence for HCpt$^2$-morphisms and the category Dir-HCpt$^2$.

\subsection{Singular homology groups are Ch-shape invariant}

As is well-known, each Hausdorff space pair $(X,X_0)$ has a relative homology group $H_n(X,X_0;G)$ for each $n\in\N$, where $G$ is the coefficient group that is abelian. Now consider a HCpt$^2$-system $X^*=\{(X_a,X_{0a}),[p_{aa'}],A\}$. Then given an abelian group $G$, similar to \cite{Mar2}, since the homotopy maps induce the same homomorphism between singular homology groups, we define
$$H_n(X^*;G)=\{H_n(X_a,X_{0a};G),p_{aa'n},A\}$$
to be a direct system in the category of groups, where if $a\leq a'$ then
$$p_{aa'n}:H_n(X_a,X_{0a};G)\ra H_n(X_{a'},X_{0a'};G)$$
is the homomorphism induced by $p_{aa'}$, as in \cite{Hat}. We often omit the group symbol $G$ in the following, writing simply $H_n(X^*)$, $H_n(X_a,X_{0a})$.

Consider the HCpt$^2$-morphism $F=([f_a],f):X^*\ra Y^*$. We can naturally obtain for each $n\in\N$, a morphism of direct systems
$$H_n(F)=(f_{an},f):H_n(X^*)\ra H_n(Y^*)$$
defined such that $f:A\ra B$ is an order-preserving map and if $a\in A$ then
$$f_{an}:H_n(X_a,X_{0a})\ra H_n(Y_{f(a)},Y_{0f(a)})$$
is induced by the map $f_a:(X_a,X_{0a})\ra(Y_{f(a)},Y_{0f(a)})$.
Thus, the identity map $1_{X^*}$ induces the identity morphism $H_n(1_{X^*})=1_{H_n(X^*)}$ for each $n\in\N$; the morphism induced by a composition is the composition of the induced morphisms, i.e., $H_n(FG)=H_n(F)H_n(G)$. This defines a covariant functor $H_n$ from the category dir-HCpt$^2$ to dir-$\sG$ for each $n\in\N$.
\bt\label{t6.1} The functor $H_n$ preserves homotopy, i.e., if $F\simeq G:X^*\ra Y^*$, then $H_n(F)\simeq H_n(G):H_n(X^*)\ra H_n(Y^*)$, $n\in\N$.\et

\bo The proof is as in Theorem 2.1 in \cite{Sand1} with the functor $\pi$ therein replaced by $H_n$. We omit the details.\eo
\bco\label{co6.1} If $X^*\simeq Y^*$, then $H_n(X^*)\simeq H_n(Y^*)$, $n\in\N$.\eco

It follows immediately from Corollary \ref{co6.1} that the functor $H_n:\mb{dir-HCpt}^2\ra\mb{dir-}\sG$ induces a functor $[H_n]:\mb{Dir-HCpt}^2\ra\mb{Dir-}\sG$ for each $n\in\N$.

Now we consider a Hausdorff space pair $(X,X_0)$. A compact Hausdorff subspace pair of $(X,X_0)$ is a Hausdorff space pair $(K,K_0)$ such that $K$, $K_0$ are compact subsets of $X$, $K_0\ss K\ss X$ and $K_0\ss X_0$. All the compact Hausdorff subspace pairs construct a directed set $c(X,X_0)$ such that $(K,K_0)\leq(K',K'_0)$ implies $(K,K_0)\ss(K',K'_0)$, i.e., $K\ss K'$ and $K_0\ss K'_0$. Associated with $(X,X_0)$ we have a HCpt$^2$-system
\be\label{6.4}C(X,X_0)=\{(K,K_0),[i_{(K,K_0)(K',K'_0)}],c(X,X_0)\},\ee
where if $(K,K_0)\leq(K',K'_0)$ then $i_{(K,K_0)(K',K'_0)}$ is the inclusion map from $(K,K_0)$ to $(K',K'_0)$. Notice that these settings are almost the same as those of homotopy groups in Subsection \ref{s4.3}. Similarly we have the following result:
\bt\label{t5.3} If ${\rm Sh}_{\rm CH}(X,X_0)={\rm Sh}_{\rm CH}(Y,Y_0)$, then $H_n(X,X_0)\approx H_n(Y,Y_0)$ for each $n\in\N$.\et

One can prove Theorem \ref{t5.3} as the process given in the proof of Theorem \ref{t4.8}. Also one can refer to this result in \cite{Hat} as follows, and Theorem \ref{t5.3} follows immediately from it.

\bp\label{p6.1}\cite{Hat} If a space $X$ is the union of a directed set of subspaces $X_a$ with the property that each compact set in $X$ is contained in some $X_a$, then the natural homomorphism $\disp\limra H_n(X_a;G)\ra H_n(X;G)$ is an isomorphism for all $n\in\N^+$ and coefficient groups $G$.\ep

\subsection{Shape cohomology groups}

Now that we have given the conclusion on the homology property of CH-shape, it is natural to consider the corresponding cohomology groups. However the cohomology groups concerning CH-shape has no such good properties.

Since we have a relative cohomology group $H^n(X,X_0)$ for each $n\in\N$ and each Hausdorff space pair $(X,X_0)$, similar to the homology case, each HCpt$^2$-system $X^*=\{(X_a,X_{0a}),[p_{aa'}],A\}$ corresponds to an inverse system in the category of groups
$$H^n(X^*)=\{H^n(X_a,X_{0a}),p_{aa'}^n,A\},$$
where $p_{aa'}^n:H^n(X_{a'},X_{0a'})\ra H^n(X_a,X_{0a})$.
Given an HCpt$^2$-morphism
$$F=([f_a],f):X^*\ra Y^*=\{Y_b,[q_{bb'}],B\},$$
we have for each $n\in\N$, an induced morphism of inverse systems
$$H^n(F)=(f_a^n,f):H^n(Y^*)\ra H^n(X^*)$$
such that $f_a$ induces $f_a^n:H^n(Y_{f(a)},Y_{0f(a)})\ra H^n(X_a,X_{0a})$ and $p^n_{aa'}f_{a'}^n=f_a^n q_{f(a)f(a')}^n$.

Define $H^n(1_{X^*})=1_{H^n(X^*)}$ to be induced by $1_{X^*}$ and the morphism induced by a composition as the composition of induced morphisms, i.e., $H^n(FG)\!=\!H^n(G)H^n(F)$. Thus one obtains a contravariant functor from dir-HCpt$^2$ to dir-$\sG$.

Similarly, for any $n\in\N$,\benu\item[(i)]the functor $H^n$ preserves homotopy; and \item[(ii)]$H^n(X^*)\simeq H^n(Y^*)$ whenever $X^*\simeq Y^*$.\eenu
\vs

Now we consider a Hausdorff space pair $(X,X_0)$ whose associated HCpt$^2$-system is $C(X,X_0)$ written as (\ref{6.4}). Then we define the {\em $n$-th CH-shape cohomology group} of $(X,X_0)$ as
$$\hat H^n(X,X_0)=\limla H^n(C(X,X_0)),$$
where $\disp\limla$ is the inverse limit, see \cite{Mar2}. Here we use the symbol \^{ } to indicate the duality to that obtained in the inverse system approach given by Marde\v si\'c and Segal \cite{Mar2}, where they use \v{ } in commemoration of \v Cech.

With the properties of inverse limits, similar to the case of CH-shape homology groups, we can define the {\em induced homomorphism} $\hat\vp^n:\hat H^n(Y,Y_0)\ra\hat H^n(X,X_0)$ of the CH-shape cohomology groups from the continuous map $\vp:(X,X_0)\ra(Y,Y_0)$, and the {\em coboundary homomorphism} $\hat\de^n$ from $\hat H^n(X_0)$ to $\hat H^n(X,X_0)$ for each $n\in\N$.

Disappointingly, the CH-shape cohomology group is not always isomorphic to the singular cohomology group for an arbitrary Hausdorff space. But when $(X,X_0)$ is a compact Hausdorff space pair, $(X,X_0)$ itself is a CS-cofinal subset of $c(X,X_0)$. Then by the equivalence of $\{(X,X_0),[1_{(X,X_0)}]\}$ and $C(X,X_0)$ in dir-HCpt$^2$, one easily sees $H^n(X,X_0)$ and $\hat H^n(X,X_0)$ are isomorphic, since isomorphisms of inverse systems induce isomorphisms of inverse limits, see \cite{Eil,Mar2}. Similarly to the case of homology groups, we have the following consequence.
\bt For a compact Hausdorff space pair, CH-shape cohomology theory and singular cohomology theory coincide.\et

For the general case, one can check Eilenberg and Steenrod's seven axioms for CH-shape cohomology groups. Using the theory of inverse systems, one can verify the axioms are perfectly satisfied trivially except the naturality axiom, excision axiom and exactness axiom. In particular, the CH-shape may not satisfy the exactness axiom generally. In the following we check the naturality and excision axioms:
\benu\item[]{\bf Naturality Axiom:} If $\vp:(X,X_0)\ra(Y,Y_0)$ is a map, then
\be\label{6.1}\hat\de^n\widehat{(\vp|_{X_0})}^{n+1}=\hat\vp^{n+1}\hat\de^n.\ee
\item[]{\bf Excision Axiom:} If $(X,X_0)$ is a Hausdorff space pair, and an open subset $W$ of $X$ satisfies $\ol{W}\ss{\rm int}X_0$, then the homomorphism
$$\hat i^n:\hat H^n(X,X_0)\ra\hat H^n(X\sm W,X_0\sm W)$$
induced by the inclusion $i:(X\sm W,X_0\sm W)\ra(X,X_0)$ is a isomorphism.\eenu

\bt The CH-shape cohomology groups satisfy the naturality axiom.\et
\bo Given $\vp$ and $\vp|_{X_0}$, we can define $F=([f_{(K,K_0)}],f):C(X,X_0)\ra C(Y,Y_0)$ and $G=([g_{K_0}],g):C(X_0)\ra C(Y_0)$. Correspondingly we have
$$H^n(F)=(f_{(K,K_0)}^n,f): H^n(C(X,X_0))\ra H^n(C(Y,Y_0))$$
$$\mb{and }\;H^n(G)=(g_{K_0}^n,g): H^n(C(X_0))\ra H^n(C(Y_0))$$
for each $n\in\N$. In the consideration of boundary homomorphism, we have the following morphisms for each $n\in\N^+$,
$$\De^n_X=(\de_{(K,K_0)}^n,\phi):H^{n}(C(X_0))\ra H^{n+1}(C(X,X_0))$$
$$\mb{and }\;\De^n_Y=(\de_{(L,L_0)}^n,\phi):H^{n}(C(Y_0))\ra H^{n+1}(C(Y,Y_0)),$$
where $\phi$ maps a pair to the second component, i.e., $\phi((K,K_0))=K_0$. To check the naturality axiom (\ref{6.1}), it suffices to check
\be\label{6.2}\De_X^n H^{n+1}(G)\approx H^{n+1}(F)\De_Y^n\hs\mb{in dir-}\sG.\ee

Indeed, given any $(K,K_0)\in c(X,X_0)$, by the naturality axiom of singular cohomology groups, we have the following commutative diagram:

\begin{picture}(0,80)
\put(150,65){\vector(-1,0){70}}
\put(150,8){\vector(-1,0){70}}
\put(43,15){\vector(0,1){43}}
\put(180,15){\vector(0,1){43}}
\put(98,70){\small$f_{(K,K_0)}^n$}
\put(108,-1){\small$g_{K_0}^n$}
\put(19,62){\small$H^n(K,K_0)$}
\put(156,62){\small$H^n(\vp(K),\vp(K_0))$}
\put(19,5){\small$H^{n+1}(K_0)$}
\put(156,5){\small$H^{n+1}(\vp(K_0))$}
\put(10,33){\small$\de_{(K,K_0)}^n$}
\put(185,33){\small$\de_{(\vp(K),\vp(K_0))}^n$}
\end{picture}

Note also $\phi f=g\phi$. (\ref{6.2}) is a trivial deduction.\eo

\bt\label{t6.1A} The CH-shape cohomology groups satisfy the excision axiom.\et

To show Theorem \ref{t6.1A}, we need  the following  result.
\bp\label{p6.2} (\cite{Eil}) Let $X_*$ and $Y_*$ be two inverse systems over directed sets $A$ and $B$ and $F=(f_b,f):X_*\ra Y_*$ be a morphism. If there is a directed subset $B'$ of $B$ such that (i) $B'$ is cofinal in $B$, (ii) $f(B)$ is cofinal in $A$, and (iii) $f_b$ is a isomorphism of $X_{f(b)}\ra Y_b$ for each $b\in B'$. Then $\disp f_{\8}=\limla F$ is an isomorphism from $X_\8$ to $Y_\8$.\ep

\noindent{\it Proof of Theorem \ref{t6.1A}.}\hs First let $A=c(X,X_0)$, $B=c(X\sm W,X_0\sm W)$ and $f:A\ra B$ such that $f((K,K_0))=(K\sm W,K_0\sm W)$. Then give any $b\in B$ a new index $b_a$ when $f(a)=b$. Thus we have a new directed set $C=B_A=\{b_a:b\in B,\,a\in A\mb{ and }f(a)=b\}$ where if $a\leq a'$ then $(f(a))_a\leq(f(a'))_{a'}$. Let $A^*_W=\{c,[q_{cc'}],C\}$ be a new HCpt$^2$-system, where if $c\leq c'$ then $q_{cc'}$ is the inclusion map.

Note that any $b_a\in C$ is actually $b$ itself with an index $a$. We define a map $g:B\ra C$ such that $g(b)=b_b$ and identities $g_b=1_b$, which form a morphism $G=([g_b],g):C(X\sm W,X_0\sm W)\ra A_W^*$. Then considering the singular homology group, we have a morphism in inv-$\sG$ for each $n\in\N$:
$$H^n(G)=(g^n_b,g):H^n(A_W^*)\ra H^n(C(X\sm W,X_0\sm W)),$$
where $g^n_b:H^n(g(b))\ra H^n(b)$ is induced by the identity $g_b$ and hence is isomorphic for all $n\in\N$ and $b\in B$.

Trivially $g(B)$ is cofinal in $C$, and one can apply Proposition \ref{p6.2} and obtain that
$$\hat g^n=\limla H^n(G):\limla H^n(A_W^*)\ra\hat H^n(X\sm W,X_0\sm W).$$
is an isomorphism.

By the construction above, one easily sees if we set $j:C\ra A$ with $j((f(a))_a)=a$ then $j$ is an order isomorphism from $C$ to $A$, and hence the inverse map $j^{-1}$ is also an order isomorphism. Moreover, we define $j_{b_a}=i_{ba}$ to be the inclusion of $b=f(a)$ into $a$. Then we obtain a morphism $J=([j_c],j):A_W^*\ra C(X,X_0)$, which induces a morphism in inv-$\sG$ for $n\in\N$,
$$H^n(J)=(j_c^n,j):H^n(C(X,X_0))\ra H^n(A_W^*).$$

Now we define a directed subset $\cA_W$ of $A$ as
$$\cA_W=\{(K,K_0)\in A:K\cap\ol{W}\ss{\rm int}_KK_0\},$$
where ${\rm int}_KK_0$ means the interior of $K_0$ with respect to $K$, i.e., there is a closed subset $U$ of $X$ with $U\cap K=K_0$ such that ${\rm int}U\cap K={\rm int}_KK_0$.

{\bf Claim:} $\cA_W$ is cofinal in $A$ and $j^{-1}(\cA_W)$ is cofinal in $C$.

Indeed, for any $(K,K_0)\in A$, since $K$ is compact and $X$ is Hausdorff, we have $K$ is normal and hence there is an open subset $U$ in $K$ such that
$$K\cap\ol{W}\ss U\ss\ol{U}\ss{\rm int}A\cap K.$$
Then set $K'_0=K_0\cup(K\cap\ol{U})$ and we have $(K,K_0)\ss(K,K'_0)\in\cA_W$, which implies the cofinality of $\cA_W$ in $A$. Moreover, given any $b_a\in C$, by the cofinality of $\cA_W$ in $A$, we have $a'\in\cA_W$ such that $a\leq a'$ and so
$$b_a\leq (f(a'))_{a'}=j^{-1}(a')\in j^{-1}(\cA_W),$$
which guarantees the cofinality of $j^{-1}(\cA_W)$ in $C$ and confirms the claim.
\vs

Given any $c=b_a\in j^{-1}(\cA_W)$, let $a=(K,K_0)$ and $b=(K\sm W,K_0\sm W)$. Since $a\in\cA_W$, by the excision property of relative singular homology groups, the induced homomorphism
$$j^n_c:H^n(K,K_0)\ra H^n(K\sm W,K_0\sm W)$$
by $j_c$ is an isomorphism for each $n\in\N$. Now by the claim and Theorem \ref{t6.1}, we have
$$\hat j^n=\limla H^n(J):\hat H^n(X,X_0)\ra\limla H^n(A_W^*)$$
is an isomorphism.

Now note $g_b$ is an identity and $j_c$ is an inclusion. It is clear that $jg:B\ra A$ is an inclusion and $j_cg_b:b\ra b$ is actually an identity, and so up to the homotopy, $JG=([j_cg_b],jg)$ is indeed the natural inclusion morphism induced by $i$. Hence we have
$$\hat i^n=\limla(H^n(JG))=\limla(H^n(G)H^n(J))=\limla H^n(G)\limla H^n(J)=\hat g^n\hat j^n$$
is an isomorphism, which assures the excision axiom.\qed

\br As a completion of shape theory (Marde\v si\'c and Segal \cite{Mar2}), CH-shape has more interesting topics (Whitehead theorem, movability, stability and so on) and applications (Morse theory, dynamical systems) worth deep studying.\er

\section*{Acknowledgement}

The author Jintao Wang wishes to express his sincere gratitude to Professor Nikita Shekutkovski for his many helpful e-mail discussions about his research, which opened for the author the door to an abundant field of shape theory and shape-like theory.


\bibliographystyle{plain}

\end{document}